\documentclass[12pt,leqno,twoside]{amsart}
\usepackage{amssymb}
\usepackage{amsmath,amsfonts}
\usepackage{amscd}
\usepackage{latexsym}
\usepackage{verbatim}
\usepackage{epsfig}

\newtheorem{theorem}{Theorem}[section]

\theoremstyle{definition}

\theoremstyle{remark}
\newtheorem{remark}[theorem]{\sc Remark}
\theoremstyle{remark}

\theoremstyle{remark}

\theoremstyle{remark}

\theoremstyle{remark}

\theoremstyle{remark}

\renewcommand{\Bbb}{\mathbb}
\newcommand{\cal}{\mathcal}

\newcommand{\bR}{{\Bbb R}}
\newcommand{\bC}{{\Bbb C}}

\newcommand{\bP}{{\Bbb P}}
\newcommand{\bN}{{\Bbb N}}
\newcommand{\bZ}{{\Bbb Z}}

\def\growarrow#1{
  \setbox1=\hbox{ $\scriptstyle #1$\ }

\mathop{\smash{\hbox to \wd1{\rightarrowfill}}
          \vphantom\rightarrow}\limits^{#1}}

\begin{document}

\title[The Levi problem on Stein spaces with singularities. A survey]
{The Levi problem on Stein spaces with singularities. A survey}

\author{Mihnea Col\c toiu}
\address{{\rm M. Col\c toiu}:
Institute of Mathematics of the Romanian Academy, P.O. Box 1-764, RO-014700, Bucure\c sti, Romania.}
\email{Mihnea.Coltoiu@imar.ro}





\maketitle

\section{A brief history of the smooth case}

In 1910  E. E. Levi \cite{Lev} noticed that a domain of holomorphy
$\Omega$ in $\bC^n$, with smooth $C^2$ boundary, should satisfy some pseudocovexity condition on the
boundary points. More precisely he showed that if $\rho$ is a $C^2$ defining function for the
boundary $\partial \Omega$ of $\Omega$ then the associated quadratic form $L \rho$ ( we shall
call it, as usual, the Levi form of $\rho$) is necessarily positive semi-definite on the holomorphic
tangent space $T_{z}(\partial \Omega):=\{w \in \bC^n \mid \Sigma _i ^n w_i {\partial \rho}  /  {\partial z_i(z)}=0 \}$ for any 
point $z \in \partial
\Omega$.

O. Blumenthal \cite{Blu} raised the important and difficult question on the
validity of the converse of this statement , i.e. if a domain $\Omega \subset \bC^n$
with smooth pseudoconvex boundary is necessarily a Stein domain. This problem, called
also the Levi problem, was open for a long time, untill 1953, when K. Oka \cite{O} solved it completely
in the affirmative (an indepedent proof of this result was also obtained by F. Norguet \cite{No} and by
H.J. Bremermann \cite{Brem}). More generally, K. Oka considered unbranched Riemann domaines $\pi :\Omega \to
 \bC^n$ (i.e. $\pi$ is locally biholomorphic) and proved that $\Omega$ is Stein iff $-\log d$ is a
 plurisubharmonic function on $\Omega$ where $d$ denotes the boundary distance on $\Omega$. Note that
 Riemann unbranched domains over $\bC^n$ appear naturally as domains of existence of families of holomorphc functions 
 defined on open subsets in $\bC^n$. Oka's result shows in particular that the Steiness of $\Omega$
 is a local property of its boundary. To be more precise we shall call a holomorphic map $p:\Omega \to X$ of complex
 manifolds (or, more generally, of complex spaces) a Stein morphism if every point $x \in X$ has
 a neighborhood $V=V(x)$ such that $p^{-1}(V)$ is Stein. For example, if we consider the inclusion map
 $i:\Omega \to X$ of an open subset $\Omega$ of $X$ , then $\Omega$ is called locally Stein (in $X$) iff
 the map $i$ is a Stein morphism, or equivalently each point $x \in \partial \Omega$ has a neighborhood
 $V=V(x)$ such that $V \cap \Omega$ is Stein. Oka's theorem can therefore be stated as follows: an open subset
 $\Omega \subset \bC^n$ is Stein if and only if $i$ is a Stein morphism, or, more generally, a
 Riemann unbranched morphism $\pi : \Omega \to \bC^n$ is Stein iff $\pi$ is a Stein morphism. If $\pi :
 \Omega \to \bC^n$ is a Stein morphism which has discrete fibers (even finite) but $\pi$ is not assumed
 to be locally biholomorphic (i.e. it is a branched Riemann domain) then $\Omega$ could be non Stein 
 (see \cite{F 2},\cite{C-D4}). 
 
 Having in mind Oka's result, H. Cartan at the important Colloque sur les fonctions de plusieurs variables
 held in Brussels in 1953 \cite{Car} raised the
 following problem: let $X$ be a Stein manifold and $\Omega \subset X$ a locally Stein open subset. Does it
 folow that $\Omega$ is itself Stein? (the Local Steiness Problem in the smooth case, the
 manifold case). A positive answer
 to this question has been given by H. Grauert and F. Doquier \cite{D-G}. Additionally,
  they solved the more general case of Riemann unbranched domains over Stein manifolds by proving :
  If $\pi : \Omega \to X$ is a Riemann unbranched domain over a Stein manifold $X$ and if $\pi$ is a
  Stein morphism, then $\Omega$ is itself Stein. Grauert's method \cite{D-G} is essentilally based on the
  following result: if $X \subset \bC^n$ is a complex closed submanifold then there exists an open
  neighborhood $V$ of $X$ ( which can be chosen to be Stein) and a holomorphic retract $\rho : V \to X$. 
  This result can be used to reduce the general case of Riemann unbranched domains over Stein
  manifolds to the case of Riemann unbranched domains spread over $\bC^n$. Such a holomorphic retract
  does not exist if $X$ has singularities. In fact Rossi \cite{Ro} showed that the existence of a
  holomorphic retract onto $X$  as above, even of  a neighborhood of $X$ minus a point, onto $X$ minus
  a point, implies the smoothness of $X$. Therefore, in order to study the Levi problem in the case
  of Stein spaces with singularities, one needs other methods and new ideas.
  
  \section {The most important questions related to the Levi problem for singular Stein spaces}
  
  The main open problem, related to the Levi problem in  the singular case (note that, due to its
importance, a new  section regarding the Levi problem in the singular case was introduced
  in the 2000 AMS classification 32C55 : The Levi problem in complex spaces; generalizations) is the Local Steiness Problem, or the
  singular Levi problem, which can be stated as follows :
  
  \smallskip 
  
   {\bf Question 1}
  
  Let $X$ be a Stein space and $D \subset X$ an open subset which is locally Stein. Does it follow that
  $D$ is itself Stein?
  
  \smallskip
  
  A more general question, in Oka's context of Riemann unbranched domains , is the following :
  
  {\bf Question 1'}
  
  Let $X$ be a Stein space and $\pi:\Omega \to X$ a Riemann unbranched domain such that $\pi$ is a
  Stein morphism. Does it follow that $\Omega$ is itself Stein?
  
  \smallskip
  
  If $\pi$ is the inclusion map, one gets the particular case of the Question 1.
  
  \smallskip
  
  Another importat open problem on Stein spaces with singularities is the Union Problem:
  
  {\bf Question 2}
  
  Let $X$ be a Stein space and $D=\bigcup_{n \in \bN} D_n$ an increasing union of Stein open sets. Does it
  follow that $D$ is itself Stein ?
  
  The answer to the Union Problem is yes if $X=\bC^n$ \cite{B-S} (one does not need necessarily
  the solution to the
  Levi problem, but the distance to the boundary is very important in the proof) and more generally
  if $X$ is a Stein manifold  using Grauert's result on the existence of the holomorphic retraction \cite{D-G}.

  \section{ The state of the art for the problems Q1, Q1', Q2 for isolated singularities}
  
  In 1964 in \cite{A-N} A. Andreotti and R. Narasimhan solved Q 1 for Stein spaces with isolated 
  singularities. Therefore they proved that a locally Stein open subset $D$ of a Stein space
  $X$ with isolated singularities is itself Stein. Their method is of "projective" type, namely they
  realize $D$ as a suitable finite union of Riemann unbranched domains $\Phi_j:D_j \to \bC^n$, $n=\dim X$, 
  $D_j \subset X$, and, using the corresponding boundary distances $d_j$ on $D_j$, they are able
  to construct, by a patching technique, a strongly plurisubharmonic continuous exhaustion
  function $\phi:D \to \bR$, and consequently, by a result due to H. Grauert \cite{G 1}, generalized
  in the singular case by R. Narasimhan \cite{Nar}, it follows that $D$ is a Stein space.
  
  This "projective" method cannot be applied to the more general case of Q1', i.e. for Riemann
  unbranched domains over Stein spaces with isolated singularities because one cannot control
  the behaviour of the constructed function in the vertical direction, i.e. in the 
  fiber direction of the map $\pi$. However it was proved by M. Col\c toiu and K. Diederich \cite{C-D4}
  that the answer to Q1' is also positive for isolated singularities, i.e. one has :
 
 \begin{theorem}
 
 Let $X$ be a Stein space with isolated singularities and  $\pi:\Omega \to X$ a Riemann unbranched domain
 such that $\pi$ is a Stein morphism. Then $\Omega$ is itself Stein.

 \end{theorem}
  
  For the proof, in order to avoid the difficulties in the construction of a function  with nice behaviour ( Levi 
  form bounded from below ) in the vertical direction, in \cite{C-D4} it is considered the pull-back 
  of the given Riemann domain on a resolution of singularities $\tau: \tilde X \to X$ and the
  existence theorem due to M. Col\c toiu and N. Mihalache \cite{C-M2} of a strongly plurisubharmonic 
  function $\phi: \tilde X \to [-\infty,\infty)$ which is $-\infty$ exactly on the exceptional
  set of the desingularization $\tilde X$ ( for the basic theory of exceptional sets
  see \cite{G 2}).  Finally one uses a classical patching technique for
  strongly plurisubharmonic functions with bounded differences.
A wrong proof of Theorem 3.1. was given by V. Vajaitu \cite{V} which is based on his lemma 3.3. which
does not hold.
  
  As for Q 2 this problem (the Union Problem) is unsolved even for Stein spaces with isolated
  singularities. The most general
   result which is known to hold is the following, due to
  M. Col\c toiu and M. Tibar \cite{C-T 2}:
  
  \begin{theorem}
  
  Let $X$ be a Stein space of dimension $2$ and let $D=\bigcup_{n \in \bN} D_n$ be an increasing 
  sequence of Stein open sets. Then $D$ satisfes the discrete Kontinuit\"atssatz (the disk property).  
  
  \end{theorem}
  
  We recall that a complex space $D$ satisfies the discrete Kontinuit\"atssatz (the disk property) if
  for any sequence of maps $\varphi_{\nu} :\bar \Delta \to D$ (where $\Delta$ is the open unit disk in $\bC$)
   which are holomorphic in $\Delta$ and continuous on $\bar \Delta$, if $\cup \varphi_{\nu}(\partial \Delta)$
   is relatively compact in $D$, then $\cup \varphi_{\nu}(\bar \Delta)$ is relatively
   compact in $D$. 
  
  The proof of the above Theorem is essentially based on the recent classification of $2$ dimensional
  normal singularities due to M. Col\c toiu and M. Tibar \cite{C-T 1}. Namely, let $(X,x_0)$ be a germ of a $2$
  dimensional normal singularity and denote by $K$ the associated singularity link, i.e. in some local
  embedding, $K$ is the intersection of the boundary of a small ball (centered in $x_0$) with $X$. If the
  fundamental group $\pi_1(K)$ is a finite group, then it is well-known that $(X,x_0)$ is a quotient
  singulaity, therefore the universal covering of $X \setminus \{x_0\}$, for small $X$, is a ball
  minus a point (this singularity is of "concave" type). The main result in \cite{C-T 1} asserts that
  if $\pi_1(K)$ is an infinite group, then the universal covering of $X \setminus \{x_0\}$, for
  small $X$, is a Stein manifold (we can say that this singularity is of
  "convex" type). The proof of this assertion is divided in two steps :
  
  \smallskip 
  
  Step 1 . It is assumed that the homology group $H_1(K,\bZ)$ is an infinite group. In this case, using a
  suitable infinite "necklace" (Nori string) and a patching technique it is constructed a Stein
  covering of $X \setminus \{x_0\}$, for small $X$.
  
  \smallskip
  
  Step 2. The general case when $\pi_1(K)$ is an infinite group is reduced to the previuos step using some
  results from the classification theory of real $3$-dimensional compact manifolds in
  order to cover finitely sheeted the link $K$ by another $3$-manifold with infinite first $\bZ$-homology
  group (which will be the link of another singularity and Step 1 can be applied to this new singularity).
  
  \smallskip
  
  In connection with the disk property (discrete Kontinuit\"atssatz) of increasing unions of Stein
  domains contained in a $2$ dimensional complex space $X$ it is important to note that an
  arbitrary increasing union of Stein manifolds ( not contained in a Stein space ) might not have the
  disk property as it was proved by J.E. Fornaess \cite{F 1}. If the discrete Kontinuit\"atssatz
  condition is replaced by the continuous Kontinuit\"atssatz (the parameter indexing the discs continuously
  is $t \in \bR$) then, obviously, from the definition, it follows that  an
  arbitrary union of Stein manifolds satisfies the continuous Kontinuit\"atssatz.
  
  \section{The state of the art for arbitrary singularities}
  
  First let us briefly recall the notion of envelope of holomorphy for a domain $D$ contained as an open subset
  in a Stein space $X$ (for arbitrary complex spaces see e.g. \cite{C-D2}). The domain $D$ is said to have an
   envelope of holomorphy , say $\tilde D$, if $\tilde D$
  is a Stein space, $D \subset \tilde D$ as an open complex subspace, and every holomorphic function on $D$ extends
  uniquely to a holomorphic function on $\tilde D$ (note that, if the envelope of holomorphy exists, then
  it is unique). If $X$ is a Stein manifold then it is well-known that $D$ has always an envelope of
  holomorphy $\tilde D$ and $\tilde D$ can be realized as an unbranched Riemann domain $\pi:\tilde D \to X$. If $X$ has
  singularities (even isolated and normal) the problem of the existence of envelopes of holomorphy is
  more difficult. In his well-known article " Remarkable pseudoconvex manifolds " \cite{G 3} H. Grauert
  constructed an example of a $3$ dimensional normal Stein sapce $X$, with an isolated singularity, and
  an open subset $D \subset X$ ( which is the complement of a hypersurface ) such that $D$ has not
  an envelope of holomorphy (Grauert's example is also discussed in detail in \cite{Su}). J. Bingener
  \cite{Bi} also constructed another example with similar properties using Nagata's counterexample to the
  Hilbert 14-th problem \cite{Nag}. However his proof is quite involved. Other counterexamples have been obtained in
  \cite{C-D1} using the hypersurface section problem (see also \cite{Col 3}) or in \cite{C-D3} using the non-separation
  of the topology of the cohomology group $H^1(D,\cal O)$.
  
  In \cite{C-D2} it is proved the following result:
  
  \begin{theorem}
  
  Let $X$ be a Stein space and $D \subset\subset X$ a locally Stein open subset. Then the following two
  conditions are equivalent :
  
  1. $D$ is Stein
  
  2. $D$ has an envelope of holomorphy
  
  \end{theorem}
  
  An analogous result is shown for increasing unions of Stein domains contained in a Stein space. The proof is
  essentially based on the following theorem due to Fornaess and Narasimhan \cite{F-N} , proved using $L^2$ estimates:
   Let $D \subset \subset X$ be a locally Stein open subset contained in a normal Stein space $X$. Then for
   every point $x_0 \in (\partial D) \cap Reg(X)$ and for every sequence of points $x_n \in D$, $x_n \to x_0$,
   there exists a holomorphic function $f \in \cal O(D)$ which is unbounded on $\{x_n\}$.
   
   Concerning the envelopes of holomorphy in normal Stein spaces of dimension $2$ it seems to the author that
   the answer to the following problem is unknown: Let $X$ be a normal Stein space of dimension $2$ and $D\subset X$
   an open subset. Is it true that $D$ has an envelope of holomorphy ?
   
   Let us recall in this context, the following question raised by H. Grauert and R. Remmert \cite{G-R}: Let $X$ be a normal
   Stein space of dimension $2$ and $D \subset X$ a domain of holomorpy. Does it follow that $D$ is itself Stein?
   Under some additional topological assumptions (e.g. $D$ is locally simply connected near $Sing(X)$) it 
   follows directly from Col\c toiu-Tibar classification of normal $2$ dimensional singularities \cite{C-T 1} that
   $D$ as above is Stein.
   
   For non-isolated singularities in dimension $3$ in \cite{C-D2} it was proved the following result: Let $X$ be a
   Stein space of dimension $3$ and $H \subset X$ a hypersurface (i.e. a closed analytic subset of 
   codimension $1$). If $D=X \setminus H$ is locally Stein then $D$ is itself Stein.
   
   It is not known if an analogous result as above holds if the condition " locally Stein" is replaced by
   "an increasing union of Stein open sets".

   There is a strong connection between the  Steiness condition of a locally Stein $D$ and the question of the separation
   of the cohomology group $H^1(D,\cal O)$. Namely in \cite{J} it is shown that an open subset $D$ of a Stein
   space is Stein if it is locally Stein and if $H^1(D,\cal O)$ is separated (an analogous result holds
   for increasing unions of Stein open sets). Therefore it is interesting to decide if the following is
   true: if $X$ is a (normal) Stein space and $H \subset X$ is a hypersurface, $D:=X \setminus H$, does it
   follow that $H^1(D,\cal O)$ is separated ?
   
   The question of the separation of $H^i(X \setminus A , \cal F)$ , $A \subset X$ closed analytic subset,
   $\cal F \in Coh(X)$ was studied in detail by Siu and Trautmann \cite{S-T 1},\cite{S-T 2}, Trautmann \cite{T} (sufficient
   conditions). However in \cite{C-D3} it was constructed a normal Stein space of dimension $3$, having only
   one singular point, and a hypersurface $H \subset X$ such that for $D=X \setminus H$ the cohomology
   group $H^1(D, \cal O)$ is not separated. It would be also interesting to consider closed analytic
   subsets $A \subset \bC^n$ and to study if the cohomology groups $H^i(X \setminus A, \cal O)$ are
   separated ( for $i=1$ the answer is yes \cite{T} ). One interesting example in this context is $A \subset \bC^6$ where
   $\dim A=3$ and $A$ is the cone over the Veronese embedding $\bP^2 \to \bP^5$ and to study the vanishing of the
   cohomology groups in degree $3$ (see also W. Barth \cite{Ba} and M. Col\c toiu \cite{Col 4}).
   In view of this discussion it would be important to decide the answer to the question : given a
   Stein space, $dim X \geq 4$, and $H \subset X$ a hypersurface ( closed analytic subset of 
   codimension $1$) such that $D:=X \setminus H$ is locally Stein, does it follow that $D$ is itself
   Stein ? 
   
   \section{The connection between Levi's problem and the hypersurface section problem}
   
   Related to the Local Steiness Problem, is the following question " The hypersurface section problem"
   considered for the first time by J. E. Fornaess and R. Narasimhan \cite{F-N} (under some
   additional cohomological vanishing assumptions) .
   
   This problem can be stated as follows :
   
   \smallskip
   
   {\bf Question}
   Let $X$ be a Stein space, $dim X \geq 3$, and $D \subset X$ an open subset such that the intersection
   $D \cap H$ is Stein for any hypersurface $H \subset X$. Does it follow that $D$ is Stein ?
   
   \smallskip
   
   A counterexample of dimension $3$, with $X$ normal, having only one singular point, and $D$ is the
   complement of a hypersurface $A \subset X$, has been constructed in \cite {C-D1} (see also \cite{Col 3}).
   An important tool in this construction are the line bundles which are topologically trivial but
   none of their power is analytically trivial (which were studied for the first time by
   H. Grauert \cite{G 3})  and a result of R. R. Simha \cite{Sim} about the Steiness of the complement
   of a curve in a $2$ dimensional normal Stein space. Later, H. Brenner \cite{Bren} obtained another 
   $3$ dimensional counterexample $X$, but with $X$ non-normal, using "forcing equations" and the
   result of R. R. Simha, with a construction which is more algebraic then geometric. He
   communicated to the author that the normalization of his counterexample works also for
   the hyperintersection problem and has only one singular point, but the computation to prove
   that the normalization has only one singular point is quite involved ( polynomials of degree $12$).
   
   The connection between the "Local Steiness Problem" LSP (The Levi Problem) and the "Hypersurface
   section problem" HSP is the following: HSP implies LSP, which follows by induction on $\dim(X)$
   since LSP is known to hold if $\dim(X)=2$ (\cite{A-N}). However, as we already remarked, HSP does
   not hold if $\dim(X)=3$. It would be interesting to construct counterexamples to HSP with
   $\dim(X) \geq 4$. This is a much more difficult problem than the $3$-dimensional case (absence of a
   Simha type result if $\dim > 2$). As remarked in \cite{Col 3} in order to construct a counterexample
   $X$ with $\dim(X) \geq 4$, it suffices to construct a compact projective algebraic space $M$, 
   $\dim(M) \geq 3$, and an open subset $U \subset M$ such that :
   
   1)  $U$ is not Stein, but the intersection $T \cap U$, of $U$ with every hypersurface $T \subset M$, 
   is Stein
   
   2) $U$ is weakly pseudoconvex, i. e. $U$ admits a smooth plurisubharmonic exhaustion function
   
   \smallskip
   
   If $X$ is a Stein manifold the answer to HSP is known to be affirmative. For $X=\bC^n$ this problem
   was solved by P. Lelong \cite{L}, and the general case of Stein manifolds follows easily from this
   case.
   
   \section{Some final remarks and conjectures}
   
   Let $\pi:X \to Y$ be a proper holomorphic map of complex spaces. We recall that $\pi$ is called relatively
   ample if there is a holomorphic line bundle $p:L \to X$ over $X$ such that the restriction of $L$ to the
   fibers of $\pi$ is an ample (positive) line bundle.
   
   In connection with the Local Steiness Problem we raise the following :
   
   {\bf Question A}
   
   Let $\pi:X \to Y$ be a proper holomorphic map which is relatively ample and assume that $Y$ is a Stein
   space. Let $W \subset X$ be an open subset such that the restriction of $\pi$ to $W$ is a Stein
   morphism. Does it follow that $W$ is Stein ?
   
   \smallskip
   
   A negative answer to this question would imply a counterexample to the Local Steiness Problem. This follows
   easily from the results of relative contraction \cite{K-S}.
   Let us remind in this context the Serre question  \cite{Se} : if $\pi :E \to B$ is a locally trivial 
   holomorphic fibration with Stein base $B$ and Stein fiber $F$, does it follow that the total space
   $E$ is itself Stein ?
   The first counterexample was obtaned by H. Skoda \cite{Sk} having as fiber $\bC^2$ (studied also
   by J.P.Demailly in \cite{Dem} and J.P. Rosay \cite{R}) and with a bounded Stein domain in $\bC^2$ as fiber by Coeur\'e and Loeb
   \cite{C-L}. It would be interesting to study if it is possible to obtain counterexamples to the Serre
   problem so that the automorphisms of the Stein fiber extend to automorphisms of some algebraic
   compactification of the fiber and the resulting projection map, after compactifying the fiber, is
   a projective morphism. Then one would get, according to the previous discussion ( question A), a
   counterexample to the Local Steiness Problem.
   
   A particular case of Question A is the Levi problem in a product :
   
   {\bf Question B}
   
   Let $Y$ be a Stein space (even a smooth Stein curve) and $M$ a projective algebraic manifold. Consider
   the canonical projection $\pi :M \times Y \to Y$ and let $W \subset M \times Y$ be an open subset
   such that the restriction of $\pi$ to $W$ is a Stein morphism. Does it follow that $W$ is Stein ?
   
   \smallskip
   
   If $\dim(M)=1$ and $Y$ is a Stein manifold the answer to Question B is affirmative \cite{Bru}, \cite{Mats}.
   
   Due to this remark, and taking also into account the fact ( already mentioned ) that the complement
   $D:=X \setminus A$ ( with $X$ Stein, $\dim(X)=3$, $A \subset X$ hypersurface ) is Stein if $D$ is 
   assumed to be locally Stein , it is natural to make the following :
   
   {\bf Conjecture 1} 
   
   The Levi Problem (i.e. the Local Steiness Problem) holds if $\dim(X)=3$.
   
   \smallskip

Question B is also related to the following problem : let $X$ be a Stein space with an isolated singularity
and let $Y$ be a Stein manifold. Denote $Z=X \times Y$ and let $U \subset Z$  be a locally Stein open 
subset. Does it follow that  $U$ is itself Stein ?

The main difficulty in the Levi problem in the singular case is the lack of a boundary distance $d$ 
such that $-\log d$ is plurisubharmonic for Stein open subsets. If $U \subset X$ is locally Stein
one can cover the boundary $\partial D$ of $D$ by Stein open sets $V_i$ ,$i \in \bN$ , such that
$V_i \cap D$ is Stein for each $i \in \bN$ , therefore there exist plurisubharmonic exhaustion 
functions $\phi_i:V_i \cap D \to \bR$. If $D$ is locally hyperconvex (i.e. it admits locally, 
negative plurisubharmonic exhaustion functions ) then it is possible (see e.g. \cite{C-M3}) to
achieve that the differences $\phi_i - \phi_j$ are bounded (by composing the given $\phi_i$
with suitable convex increasing functions), and consequently, by a simple patching technique, 
one gets a strongly plurisubharmonic exhaustion $\phi:D \to \bR$ ( which implies that $D$ is
Stein, by the results of H. Grauert \cite{G 1} and R. Narasimhan \cite{Nar} ). If $D$ is not locally
hyperconvex, and it is assumed only locally Stein, it seems that it is not possible to get the
plurisubharmonic local functions $\phi_i$ with bounded differences ( if we compose them with
non-convex increasing functions then, locally, their Levi form loses one positive eigenvalue, and in the
patching process one gets only $2$-completeness with corners, not even $2$-completeness).

\begin{remark}
If $D \subset X$ is locally Stein and $X$ is a Stein space , it follows by a bumping technique 
(see \cite{C-M1} ) that the cohomology groups $H^i(X,\cal O)$ vanish if $i \geq 2$. If, additionally, 
it is assumed that $H^1(D, \cal O)=0$ then one gets immedialtely, by using the Koszul complex, 
that $D$ is Stein.
Similarly it is known ( see \cite{Mar}, \cite{Sil} ) that an arbitrary increasing union of 
Stein spaces $\{X_n\}_{n \in \bN}$
is itself Stein if $H^1(X,\cal O)=0$ ( in fact it suffices to assume that $H^1(X,\cal O)$ is separated.

\end{remark}

\begin{remark}
 
By using a bumping argument one can easily see that another equivalent statement to L.S.P. (of Oka's
glueing lemma type) is the following : Let $X$ be a Stein space and $D \subset X$ an open subset. Assume
that there exists $f \in \cal O(X)$ and real constants $a < b$ such that $D \cap \{Re f < b\}$ and 
$D \cap \{Re f > a \}$ are Stein. Does it follow that $D$ is itself Stein ? ( see also \cite{A-N}
concerning Oka's Heftungslemma).

\end{remark}

Concerning the Union Problem  we already mentioned \cite{C-T 2} that for any Stein space $X$, with $dim(X)=2$, 
and for any open subset $U \subset X$ such that $U=\bigcup_{n \in \bN}U_n$, $U_n \subset U_{n+1}$, $U_n$ is
Stein for every $n \in \bN$, if follows that $U$ satisfies the discrete Kontinuit\"atssatz ( the disk
property). On the other hand any such increasing union $U \subset X$, with $X$ Stein, normal, $2$ dimensional, 
is a domain of holomorphy, therefore if the Grauert and Remmert problem \cite{G-R}, already mentioned, has a
positive answer, then it would follow that $U$ is Stein.

Therefore it is natural to make the following :

{\bf Conjecture 2}

The Union Prolem holds if $\dim(X)=2$.

\smallskip

Let us make also some remarks concerning the Union Problem for Stein spaces with isolated singularities, 
which is an open question. Suppose that $M$ is a projective algebraic manifold and $U=\bigcup_{n \in \bN} U_n$
is an increasing union of Stein open sets. By considering a negative line bundle over $M$ , and denoting
$X$ the Stein space (in fact affine algebraic) obtained by  contractiong  the null section to a 
point, we 
easily see, using some arguments involving $\bC^{*}$ fibrations, that if the Union Problem has a positive answer
for isolated singularities, then necessarily $U$ is itself Stein (i. e. the Union Problem holds for
projective algebraic manifolds). For example it would be interesting to consider the case when $U$ is
the complement of a divsor $A$ and $A$ is the limit ( in the Hausorff sense)  of a sequence of divisors 
$A_n$ whose complements are Stein for every $n \in \bN$. Does it follow that $U$ is itself Stein ?

For other problems and discussions concerning the singular Levi problem the reader is advised to consult
the papers of M. Col\c toiu \cite{Col 1} , \cite {Col 2}, J. E. Fornaess and R. Narasimhan
\cite{F-N} and the survey of Y-T. Siu \cite{Siu}.


\end{document}